\documentclass[11pt,english]{paper}
\usepackage{lmodern}

\usepackage[T1]{fontenc}
\usepackage[latin9]{inputenc}
\usepackage[letterpaper]{geometry}
\usepackage{amsthm}
\usepackage{amsmath}
\usepackage{amssymb}
\usepackage{graphicx}

\makeatletter
\numberwithin{equation}{section}
\numberwithin{figure}{section}
\numberwithin{table}{section}
\newcommand{\lyxaddress}[1]{
\par {\raggedright #1
\vspace{1.4em}
\noindent\par}
}
\theoremstyle{plain}
\newtheorem{thm}{\protect\theoremname}[section]

\@ifundefined{showcaptionsetup}{}{%
 \PassOptionsToPackage{caption=false}{subfig}}
\usepackage{subfig}
\makeatother

\usepackage{babel}
\providecommand{\theoremname}{Theorem}

\begin{document}

\title{Accuracy and stability of inversion of power series}

\author{Raymundo Navarrete and Divakar Viswanath}

\maketitle

\lyxaddress{Department of Mathematics, University of Michigan (raymundo/divakar@umich.edu).}
\begin{abstract}
This article considers the numerical inversion of the power series
$p(x)=1+b_{1}x+b_{2}x^{2}+\cdots$ to compute the inverse series $q(x)$
satisfying $p(x)q(x)=1$. Numerical inversion is a special case of
triangular back-substitution, which has been known for its beguiling
numerical stability since the classic work of Wilkinson (1961). We
prove the numerical stability of inversion of power series and obtain
bounds on numerical error. A range of examples show these bounds to
be quite good. When $p(x)$ is a polynomial and $x=a$ is a root with
$p(a)=0$, we show that root deflation via the simple division $p(x)/(x-a)$
can trigger instabilities relevant to polynomial root finding and
computation of finite-difference weights. When $p(x)$ is a polynomial,
the accuracy of the computed inverse $q(x)$ is connected to the pseudozeros
of $p(x)$.
\end{abstract}

\section{Introduction}

Suppose $p(x)$ is the power series $1+b_{1}x+b_{2}x^{2}+\cdots$.
We consider the numerical accuracy and stability of computing $q(x)=1+c_{1}x+c_{2}x^{2}+\cdots$
such that $q(x)=1/p(x)$. No assumption is made regarding the convergence
of either series. It is only required that the Cauchy product $p(x)q(x)=1$.

The algorithm for inverting power series is especially simple. It
is a specialized form of triangular back substitution. To find $c_{k}$
we use $c_{k}=-b_{k}-\sum_{j=1}^{k-1}c_{j}b_{k-j}$ in the order $k=1,2,3,\ldots$ 

Inversion of power series arises as an auxiliary step in polynomial
algebra, Hermite interpolation and computations related to Padé approximation
\cite{ButcherCorless2011,SadiqViswanath2013}, where a knowledge of
its numerical properties would be useful. Yet the algorithm itself
is so simple that it appears appropriate to state that the numerical
properties of triangular back substitution are especially subtle.
In his classic paper \cite{Wilkinson1961}, Wilkinson provided a rounding
error analysis of triangular back-substitution and remarked that the
algorithm itself appeared more accurate than the error bounds. In
particular, while the bounds predict relative error proportional to
the condition number, the actual errors appear independent of condition
numbers. Higham \cite{Higham1998,Higham2002} has refined and extended
Wilkinson's analysis.

It is well-known that some natural and obvious methods for basic tasks
such as computing the standard deviation or solving a quadratic equation
are numerically unstable \cite{Higham2002}. In Section 2.1, we begin
with a similar phenomenon that arises when a polynomial $p(x)$, for
which $x=a$ is a root satisfying $p(a)=0$, is deflated to compute
$q(x)=p(x)/(x-a)$. This step arises in polynomial root finding as
well as the computation of finite difference weights. We show that
an obvious method for deflating by a root has a catastrophic numerical
instability. Indeed a general method for calculating spectral differentiation
matrices, implemented by Weideman and Reddy \cite{WeidemanReddy2000},
suffers from this instability as the order of the derivative increases,
as shown earlier in \cite{SadiqViswanath2014}. Here, in Section 2.1,
we show why the instability arises in a seemingly harmless situation.

The problem of deflating by a root is related to but not exactly the
same as that of inverting a power series. In Section 2.2, we consider
the special case of inverting a quadratic. These two problems of Section
2 bring to light some of the issues that arise in inverting power
series in a relatively transparent manner.

The notion of pseudozeros due to Mosier \cite{Mosier1986} (who called
them root neighborhoods) and in greater generality to Toh and Trefethen
\cite{TohTrefethen1994} may be invoked to shed further light on rounding
errors that arise during inversion of polynomials. The rounding errors
in coefficients of the inverse series are eventually dominated by
the polynomial root closest to the origin. However, the bounds based
on pseudozeros and condition numbers are not good. Much like algorithms
for polynomial root finding, condition numbers are derived thinking
one root at a time. In contrast, the perturbative errors in the roots
are finely correlated and the correlation in errors leads to much
better accuracy than the bounds indicate.

In Section 3, we give better bounds for the rounding errors that arise
while inverting power series. These bounds imply the numerical stability
of power series inversion. Computations that utilize extended precision
arithmetic (with $100$ digits of precision) show that the bounds
are quite good. There is no significant gap between numerical condition
and actual errors unlike the situation with triangular matrices.

A significant contribution to explain the puzzle raised by Wilkinson
\cite[p. 320]{Wilkinson1961}, namely the observed independence of
relative errors from condition numbers in triangular back substitution,
was made by Stewart \cite{Stewart1997}. Stewart has noted that triangular
matrices that arise from Gaussian elimination or QR factorization
are likely to be rank-revealing (in a sense explained in Section 3).
For such matrices, Stewart has proved that the ill-conditioning can
be eliminated using row scaling, thus partially explaining Wilkinson's
observation. The triangular Toeplitz matrices associated with power
series are typically not rank-revealing but can be so in some situations,
as shown in Section 3, but in these situations power series inversion
is well-conditioned. Thus bounds for power series inversion are generally
quite good, unlike the situation with triangular matrices.

\section{Inversion of polynomials}

In this section, we first consider deflating a polynomial $p(x)$
by factoring out $(x-a)$ where $a$ is a root satisfying $p(a)=0$.
Next we look at the inversion of a quadratic polynomial and the theory
of pseudozeros. 

Following Higham \cite{Higham2002}, but with some modifications,
we set down the basic properties of floating point arithmetic. The
floating point axiom is $\mathrm{fl(x.op.y)=(x.op.y)(1+\delta)}$
where $|\delta|\leq u$. We may also write 
\[
\mathrm{fl(x.op.y)=(x.op.y)/(1+\delta),}
\]
where again $|\delta|\leq u$. Here $u$ is the unit roundoff ($u=2^{-53}$
for double precision arithmetic) and $\text{op}$ may be addition,
subtraction, division, or multiplication. 

To handle the accumulation of relative error through a succession
of operations, it is helpful to introduce $\theta_{n}$ which is any
quantity that satisfies
\[
1+\theta_{n}=(1+\delta_{1})^{\rho_{1}}(1+\delta_{2})^{\rho_{2}}\ldots(1+\delta_{n})^{\rho_{n}}
\]
for $|\delta_{i}|\leq u$ and with each $\rho_{i}$ being $+1$, $-1$,
or $0$. In our usage, the $\theta$ variables are local to each usage.
So for example, if $\theta_{3}$ occurs in two different equations
or in two different places in the same equation, it is not the same
$\theta_{3}$, but each $\theta_{3}$ is a possibly different relative
error equal to the relative error from three (or fewer) operations.
If $a$ and $b$ are of the same sign, we may write $a(1+\theta_{n})+b(1+\theta_{n})=(a+b)(1+\theta_{n})$,
but not if they are of opposite signs.

It may be shown (see \cite{Higham2002}) that $|\theta_{n}|\leq\gamma_{n}$,
where $\gamma_{n}=nu/(1-nu)$, if $nu<1$. Unlike $\theta_{n}$, $\gamma_{n}$
stands for the same quantity in every occurrence. Whenever $\gamma_{n}$
is used, the assumption $nu<1$ is made implicitly. Another useful
bound is $\left(1+\gamma_{k}\right)\left(1+\gamma_{l}\right)\leq1+\gamma_{k+l}$.

\subsection{Deflation by $x-a$}

Let $p(x)=x^{n}+b_{n-1}x^{n-1}+\cdots+b_{0}$ and $p(a)=0$. Consider 

\[
\frac{x^{n}+b_{n-1}x^{n-1}+\cdots+b_{1}x+b_{0}}{x-a}=x^{n-1}+c_{n-2}x^{n-2}+\cdots+c_{1}x+c_{0}.
\]
Equating coefficients, we get the equations
\begin{align}
-ac_{0} & =b_{0}\nonumber \\
-ac_{1}+c_{0} & =b_{1}\nonumber \\
 & \vdots\nonumber \\
-ac_{n-2}+c_{n-3} & =b_{n-2}\nonumber \\
-a+c_{n-2} & =b_{n-1}.\label{eq:secn2.1-recurrence-for-ci}
\end{align}
We consider the accumulation of rounding error when these equations
are solved for $c_{i}$ in the order $c_{0},c_{1},\ldots,c_{n-2}$
using $c_{k}=\left(c_{k-1}-b_{k}\right)/a$ for $k>1$. If $\hat{c}_{k}$
is the computed quantity in floating point arithmetic, we assume inductively
that 
\[
\hat{c}_{k-1}=-\frac{b_{k-1}}{a}\left(1+\theta_{2}\right)-\frac{b_{k-2}}{a^{2}}\left(1+\theta_{4}\right)-\cdots-\frac{b_{1}}{a^{k-1}}\left(1+\theta_{2k-2}\right)-\frac{b_{0}}{a^{k}}\left(1+\theta'_{2k-2}\right).
\]
Since the recurrence $c_{k}=\left(c_{k-1}-b_{k}\right)/a$ involves
two operations, we have
\[
\hat{c}_{k}=-\frac{b_{k}}{a}\left(1+\theta_{2}\right)-\frac{b_{k-1}}{a^{2}}\left(1+\theta_{4}\right)-\cdots-\frac{b_{1}}{a^{k}}\left(1+\theta_{2k}\right)-\frac{b_{0}}{a^{k+1}}\left(1+\theta'_{2k}\right).
\]
Using $c_{k}=-\sum_{j=0}^{k}b_{j}/a^{k+1-j}$, we have the following
bound for the rounding error in $c_{k}$.
\begin{thm}
If the equations (\ref{eq:secn2.1-recurrence-for-ci}) are solved
for $c_{i}$ in the order $c_{0},c_{1},\ldots,c_{n-2}$, and $\hat{c}_{k}$
is the computed value of $c_{k}$ in floating point arithmetic, we
have 
\[
|\hat{c}_{k}-c_{k}|\leq\gamma_{2}\Biggl|\frac{b_{k}}{a}\Biggr|+\cdots+\gamma_{2k}\Biggl|\frac{b_{1}}{a^{k}}\Biggr|+\gamma_{2k}\Biggl|\frac{b_{0}}{a^{k+1}}\Biggr|\leq\gamma_{2k}\sum_{j=0}^{k}\Biggl|\frac{b_{j}}{a^{k+1-j}}\Biggr|.
\]
\label{thm:secn2.1-boundc0tocn}
\end{thm}
\begin{figure}
\subfloat[]{

\centering{}\includegraphics[scale=0.4]{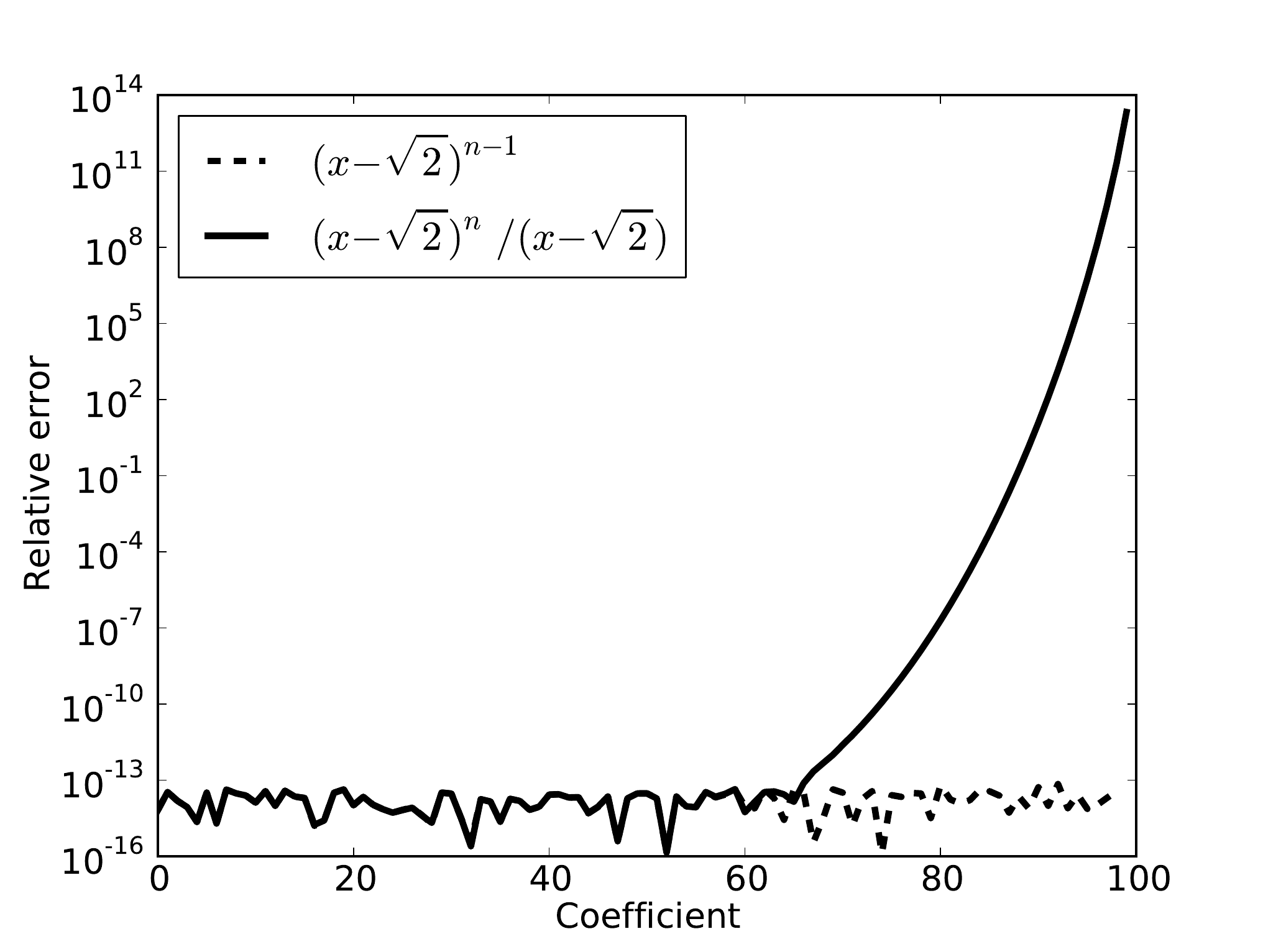}}\subfloat[]{

\begin{centering}
\includegraphics[scale=0.4]{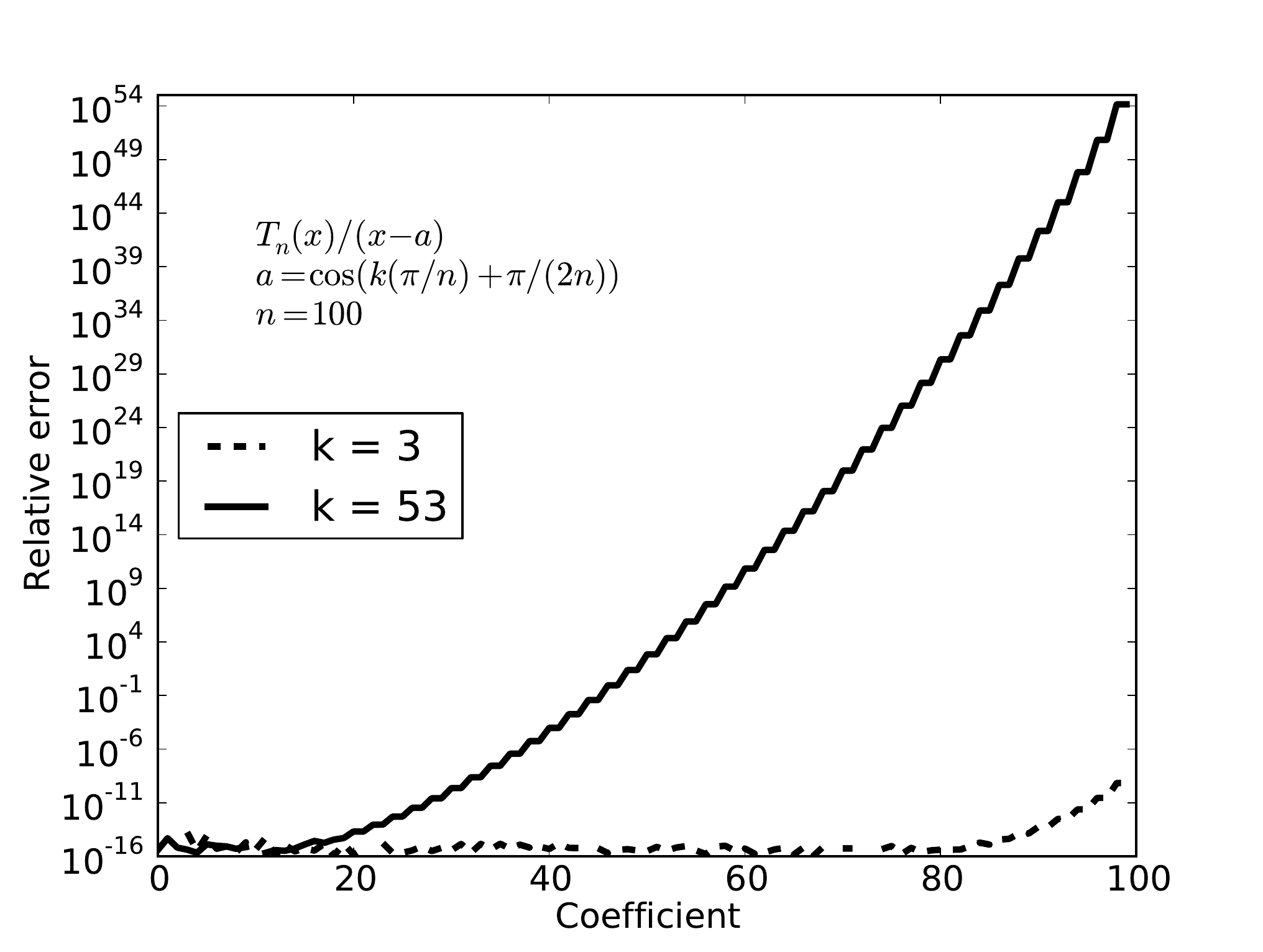}
\par\end{centering}

}

\caption{Accumulation of rounding error in the coefficients when the polynomial
$p(x)$ is deflated by $(x-a)$, $a$ being a root.\label{fig:secn2.1-error-deflation}}

\end{figure}

Within the error bound of Theorem \ref{thm:secn2.1-boundc0tocn},
there are two different mechanisms for large rounding errors. These
two mechanisms are illustrated in Figure \ref{fig:secn2.1-error-deflation}.
Figure \ref{fig:secn2.1-error-deflation}a shows the relative errors
in the coefficients of $(x-\sqrt{2})^{99}$ computed in two ways.
The first computation begins with $(x-\sqrt{2})^{100}$ and then divides
by $x-\sqrt{2}$. In the second computation, $99$ factors $x-\sqrt{2}$
are multiplied. In the first computation, it is seen that the relative
errors are initially small but begin to explode after the half way
mark. In contrast, the relative errors remain small throughout in
the second computation. 

The error bound in Theorem \ref{thm:secn2.1-boundc0tocn} corresponds
to the exact formula $c_{k}=-\sum_{j=0}^{k}b_{j}/a^{k+1-j}$. The
relative error in $c_{k}$ will be large if some of the terms of this
sum are much larger than $c_{k}$. In the binomial expansion, the
coefficients at the edges are much smaller than the ones in the middle.
Thus deflation, using the method of Theorem \ref{thm:secn2.1-boundc0tocn},
leads to large errors once we get past the middle. This is the first
mechanism for large rounding errors.

The Chebyshev polynomial $T_{n}(x)$ is defined as $\cos(n\arccos x)$
for $x\in[-1,1]$. All its $n$ roots are in the interval $[-1,1]$.
Figure \ref{fig:secn2.1-error-deflation}b shows the errors in the
coefficients of $T_{n}(x)/(x-a)$, where $a$ is a root close to $0$
and when $a$ is a root close to $1$. The errors grow explosively
for $a\approx0$ ($k=53$ in the plot), but are quite mild when $a\approx1$
($k=3$ in the plot). Here too, as indicated by Theorem \ref{fig:secn2.1-error-deflation},
there must be cancellations between the terms of $-\sum_{j=0}^{k}b_{j}/a^{k+1-j}$
for large relative errors. The cancellations can be particularly severe
when $a$ is small. This is the second mechanism for large rounding
errors.

One of the methods for computing spectral differentiation matrices
\cite{WeidemanReddy2000,Welfert1997} suffers from an instability
related to the second mechanism. This instability has been completely
fixed \cite{SadiqViswanath2014}, yet we explain exactly how it comes
about. In the original formulation \cite{WeidemanReddy2000,Welfert1997},
the connection to polynomials and root deflation is not transparent.

Equation (7) of \cite{Welfert1997}, which is the heart of the algorithm
in that paper, reads as follows:

\[
\left(D_{p+1}\right)_{k,j}=\frac{p+1}{x_{k}-x_{j}}\left(\frac{c_{k}}{c_{j}}\left(D_{p}\right)_{k,k}-\left(D_{p}\right)_{k,j}\right).
\]
Here $\left(D_{p}\right)_{k,j}$ denotes the coefficient at $x_{j}$
of the $p$-th derivative at $x_{k}$. The grid points $x_{0},x_{1},\ldots,x_{n}$
are assumed to be distinct. The $c_{i}$ are normalizing constants
extraneous to the discussion here and will be ignored. 

Let $l_{j}(x)$ denote the Lagrange cardinal function which is equal
to $1$ at $x_{j}$ and $0$ at the other grid points. The coefficients
of $u^{p}$ and $u^{p+1}$ in the polynomial $l_{j}(u+x_{k})$ 
\[
\frac{\left(u+x_{k}-x_{1}\right)\left(u+x_{k}-x_{2}\right)\cdots\left(u+x_{k}-x_{n}\right)}{u+x_{k}-x_{j}},
\]
multiplied by normalizing constants which we ignore, are equal to
the $p$-th and $(p+1)$-st derivatives of the Lagrange cardinal function
$l_{j}(x)$ evaluated at $x_{k}$, respectively (see \cite{SadiqViswanath2014}).
Similarly, the coefficient of $u^{p+1}$ of $\prod_{i=1}^{n}\left(u+x_{k}-x_{i}\right)$,
multiplied by a normalizing constant which we ignore, is equal to
the $p$-th derivative of the Lagrange cardinal function $l_{k}(x)$
evaluated at $x_{k}$\cite{SadiqViswanath2014}. Finite difference
weights are nothing but the coefficients of Lagrange cardinal functions,
suitably normalized. It follows that equation (7) of \cite{Welfert1997}
is using exactly the same recurrence as in Theorem \ref{thm:secn2.1-boundc0tocn}
and is therefore susceptible to the instability exhibited above.

Root deflation is a part of polynomial root finding algorithms such
as Jenkins-Traub \cite{JenkinsTraub1970}. In these applications,
the equations (\ref{eq:secn2.1-recurrence-for-ci}) are solved for
$c_{i}$ in the order $c_{n-1},c_{n-2},\ldots,c_{0}$. The bound in
the following theorem is proved in much the same way as the bound
in Theorem \ref{thm:secn2.1-boundc0tocn}.
\begin{thm}
If the equations (\ref{eq:secn2.1-recurrence-for-ci}) are solved
for $c_{i}$ in the order $c_{n-2},c_{n-1},\ldots,c_{0}$ (bottom
to top in (\ref{eq:secn2.1-recurrence-for-ci})), the error in the
computed quantity $\hat{c}_{k}$ satisfies the bound
\[
\big|\hat{c}_{k}-c_{k}\bigr|\leq|b_{k+1}|\gamma_{1}+|ab_{k+2}|\gamma_{3}+\cdots+|a|^{n-k-1}\gamma_{2n-2k+1}.
\]
\label{thm:secn2.1-fromcntoc0}
\end{thm}
The computation in Theorem \ref{thm:secn2.1-fromcntoc0} corresponds
to the formula $c_{k}=a^{n-k-1}+\sum_{j=k+1}^{n-1}b_{j}a^{j-k-1}$.
This appears a safer method because it is not vulnerable to the second
mechanism when $a\approx0$, and if the coefficients are well-scaled
we may assume that the roots $a$ are not too large. However, it is
still vulnerable to the first mechanism. For example, if this algorithm
is applied to deflate a factor of $\left(x-a\right)^{n}$, large errors
in the coefficients will occur for powers lower than $x^{n/2}$.

In the computation of finite difference weights, both instability
mechanism are avoided by the method of partial products \cite{SadiqViswanath2014}.
In that method the operation of deflating a polynomial by a factor
is not employed. Going by analogy, it is natural to make the suggestion
that polynomial root finding algorithms that avoid root deflation
may be more accurate for each individual root. The operation count
may be higher, but the polynomial root finding problems are puny compared
to the power of modern computers. Thus accuracy is of greater consequence.

\subsection{Inversion of a quadratic}

In a quadratic $ay^{2}+by+c$ with $ac\neq0$, we may make the change
of variables $x=sy$ and choose the scale factor $s$ to make the
coefficients of $x^{2}$ and $x$ equal in magnitude. If the coefficient
of $x^{2}$ is factored out we are left with a quadratic of the form
$x^{2}+bx\pm1$. The operations of factoring out the leading coefficient
and rescaling the variable induce minimal relative error in the computed
coefficients. Therefore as far as the accumulation of error in the
coefficients of the inverse is concerned, we are left with only two
cases:

\[
\frac{1}{x^{2}+bx\pm1}=\pm1+c_{1}x+c_{2}x^{2}+\ldots
\]
In the $-1$ case, we have $c_{1}=-b$, $c_{2}=bc_{1}-1$, and $c_{n+1}=bc_{n}+c_{n-1}$.
It follows that $c_{i}$ has the opposite sign to $b$ if $i$ is
odd and is negative if $i$ is even. There are no cancellations and
all coefficients are computed with excellent relative accuracy. Both
roots of the quadratic equation $x^{2}+bx-1=0$ are real.

The other case is with $+1$. In this case, we have 
\begin{eqnarray*}
c_{1} & = & -b\\
c_{2} & = & b^{2}-1\\
c_{3} & = & -b^{3}+2b\\
 & \vdots
\end{eqnarray*}
In general, $c_{n+1}=-bc_{n}-c_{n-1}$. Each $c_{n}$ is a polynomial
in $b$: $c_{n}=F_{n}(b)$ where $F_{n}$ is a polynomial of degree
$n$. If $\alpha$ and $\beta$ are the two distinct roots of $x^{2}+bx+1=0$
it follows that 
\begin{equation}
c_{n}=F_{n}(b)=\frac{1}{\beta-\alpha}\left(\frac{1}{\beta^{n+1}}-\frac{1}{\alpha^{n+1}}\right).\label{eq:secn2.2-cn-alpha-beta-formula}
\end{equation}
To keep the discussion simple we omit the cases $b=\pm2$ with repeated
roots. The polynomials $F_{n}$ are a version of Fibonacci polynomials
\cite{HoggattBicknell1973}. An easy induction argument using the
recurrence $c_{n+1}=-bc_{n}-c_{n-1}$ proves that the polynomial $c_{2n}=F_{2n}(b)$
has only even degree terms and that the coefficients alternate in
sign beginning with $b^{2n}$. Similarly, $c_{2n+1}=F_{2n+1}(b)$
has only odd degree terms and the coefficients alternate in sign beginning
with $-b^{2n+1}$.

If we write $c_{n-1}=F_{n-1}(b)=\sum_{k=0}^{n-1}C_{n-1,k}b^{k}$,
we may inductively assume that the computed quantity $\hat{c}_{n-1}$
is given by $\sum_{k=0}^{n-1}C_{n-1,k}b^{k}(1+\theta_{2n-2})$. Likewise,
we may inductively assume that $\hat{c}_{n}=\sum_{k=0}^{n}C_{n,k}b^{k}(1+\theta_{2n})$.
The recurrence $c_{n+1}=-bc_{n}-c_{n-1}$ implies that $C_{n+1,k}b^{k}=-\left(C_{n,k-1}b^{k-1}\right)b-C_{n-1,k}b^{k}$,
where crucially $C_{n,k-1}$ and $C_{n-1,k}$ have the same sign,
thanks to the pattern in the signs of the coefficients of $F_{n}(b)$
and $F_{n-1}(b)$. Therefore we may infer that $\hat{c}_{n+1}=\sum_{k=0}^{n+1}C_{n+1,k}b^{k}(1+\theta_{2n+2})$,
completing the induction.

The error bound 
\[
\frac{|c_{n}-\hat{c}_{n}|}{|c_{n}|}\leq\frac{|F_{n}|\left(|b|\right)}{\Bigl|F_{n}(b)\bigr|}\gamma_{2n},
\]
where $|F_{n}|$ is the polynomial with all coefficients of $F_{n}$
replaced by their absolute values, follows immediately. If we go back
to formula (\ref{eq:secn2.2-cn-alpha-beta-formula}) for $c_{n}$,
we get a sense of when the relative errors in the computed coefficients
may be large. If $|b|<2$, both roots $\alpha$ and $\beta$ of $x^{2}+bx+1$
are complex of magnitude $1$ and conjugates of each other. For certain
values of $n$, the arguments of $\alpha^{n+1}$ and $\beta^{n+1}$
will differ very nearly by a multiple of $2\pi$ and formula (\ref{eq:secn2.2-cn-alpha-beta-formula})
implies a cancellation making $F_{n}(b)$ much smaller in magnitude
than $|F_{n}|\left(|b|\right)$. The corresponding coefficients $c_{n}$
will have large relative errors.

\subsection{Connection to pseudozeros }

Let $p=p_{0}+p_{1}z+p_{n-1}z^{n-1}+z^{n}$ be a monic polynomial and
let $Z(p)=\{a_{1},a_{2},\dots,a_{n}\}$ be the set of roots of $p$.
We assume $p_{0}\neq0$. We shall connect the errors in computing
the inverse series $q(z)=1/p(z)$ to the pseudozeros of $p(z)$. The
analysis here is of conditioning not of rounding errors. We consider
another monic polynomial $\hat{p}$ close to $p$ and bound the errors
in $\hat{q}=1/\hat{p}$ using the pseudozero sets of $p$. The subscripted
variable $p_{i}$ denotes the coefficient of $z^{i}$ in $p(z)$.
Similarly $q_{i}$ denotes the coefficient of $z^{i}$ in $q(z)$. 

Pseudozero sets have been defined using the infinity norm \cite{Mosier1986}
or more general norms \cite{TohTrefethen1994}. Here we define pseudozero
sets using the maximum coefficient-wise relative error. Our definition
is close to that of \cite{Mosier1986}. Let 
\[
e(\hat{p}):=\max_{i,p_{i}\neq0}\frac{|p_{i}-\hat{p}_{i}|}{|p_{i}|}
\]
be the maximum coefficient-wise error in $\hat{p}$ relative to $p$.
The $\epsilon$-pseudozero set of $p$ in the complex plane is given
by 
\[
Z_{\epsilon}(p):=\{z\in\mathbb{C}\;:\; z\in Z(\hat{p}),\; e(\hat{p})\leq\epsilon\}.
\]
An argument in \cite{Mosier1986} (also see \cite{TohTrefethen1994})
implies that 
\[
Z_{\epsilon}(p)=\left\{ z\in\mathbb{C}\;:\;\frac{|p(z)|}{|p|(|z|)}\leq\epsilon\right\} 
\]
where $|p|$ is the polynomial with all coefficients of $p$ replaced
by their absolute values. 

Suppose $\hat{a}\in Z_{\epsilon}(p)$ and let $a\in Z(p)$, with $a=a_{i}$
for some $i$, be the root closest to $\hat{a}$. All the roots $a_{i}$
of $p(x)=0$ are assumed to be distinct, to avoid technicalities of
no value for the discussion here. Then 
\[
|a-\hat{a}|^{n}\leq\prod_{i,p(a_{i})=0}|\hat{a}-a_{i}|=|p(\hat{a})|\leq\epsilon|p|(|\hat{a}|)
\]
since $p$ is a monic polynomial. We have 
\[
|a-\hat{a}|\leq\sqrt[n]{\epsilon|p|(|\hat{a}|)},
\]
but this bound on the error is highly pessimistic. This bound is reasonably
good only if $|\hat{a}-a|\approx|\hat{a}-a_{i}|$ for every $i$,
which is very seldom the case.

Condition numbers of polynomials roots \cite{Gautschi1984,TohTrefethen1994}
may be used to derive better and less pessimistic bounds. If $a_{j}$
is a simple root of $p$ we may define
\[
\kappa(a_{j},p):=\lim_{e(\hat{p})\to0}\sup_{\hat{p}}\frac{|a_{j}-\hat{a}_{j}|}{e(\hat{p})}
\]
where $\hat{a}_{j}$ is the root of $\hat{p}$ corresponding to $a_{j}$
and $e(\hat{p})$ is the maximum relative coefficient-wise distance
of $\hat{p}$ from $p$ defined earlier. If $e(\hat{p})<\epsilon$
and $\epsilon\to0$, we have, 
\[
\frac{p(\hat{a}_{j})}{p'(a_{j})(\hat{a}_{j}-a_{j})}=\frac{(\hat{a}_{j}-a_{1})\cdots(\hat{a}_{j}-a_{j-1})(\hat{a}_{j}-a_{j+1})\cdots(\hat{a}_{j}-a_{n})}{(a_{j}-a_{1})\cdots(a_{j}-a_{j-1})(a_{j}-a_{j+1})\cdots(a_{j}-a_{n})}\to1
\]
implying $\hat{a}_{j}-a_{j}\approx p(\hat{a}_{j})/p'(a_{j})$. Therefore,
we have
\[
\kappa(a_{j},p)=\lim_{\epsilon\to0}\sup_{\hat{p},e(\hat{p})\leq\epsilon}\frac{|p(\hat{a}_{j})|/|p'(a_{j})|}{e(\hat{p}_{j})}=\frac{|p|(|a_{j}|)}{|p'(a_{j})|}
\]
noting that the inequality $|p(\hat{a})|\leq\epsilon|p|(|\hat{a}|)$
is sharp for some polynomial $\hat{p}$ with $e(\hat{p})=\epsilon$
(see \cite{Mosier1986}).

If $p$ has only distinct roots as assumed, we have
\[
q(z)=\frac{\text{Res}(q,a_{1})}{(z-a_{1})}+\cdots+\frac{\text{Res}(q,a_{n})}{(z-a_{n})}
\]
where the residue of $q$ at one of its simple poles $a_{j}$ is given
by $\text{Res}(q,a_{j})=1/[(a_{j}-a_{1})\cdots(a_{j}-a_{j-1})(a_{j}-a_{j+1})(a_{j}-a_{n})]$.
We may expand $q$ as 
\[
q(z)=\sum_{j=1}^{n}\text{Res}(q,a_{j})\left(\frac{-1}{a_{j}}\right)\sum_{k=0}^{\infty}\left(\frac{z}{a_{j}}\right)^{k}=\sum_{k=0}^{\infty}\left(\sum_{j=1}^{n}\frac{-\text{Res}(q,a_{j})}{a_{j}^{k+1}}\right)z^{k}.
\]
Let $\hat{q}=1/\hat{p}$, where $e(\hat{p})\leq\epsilon$, and let
$Z(\hat{p})=\{\hat{a}_{1},\dots,\hat{a}_{n}\}$ with $\hat{a}_{i}$
corresponding to $a_{i}$, with $\epsilon$ assumed small enough that
the correspondence may be set up. The error in the coefficient of
$z^{k}$ is
\[
(q-\hat{q})_{k}=\sum_{j=1}^{n}\left(\frac{\text{Res}(\hat{q},\hat{a}_{j})}{\hat{a}_{j}^{k+1}}-\frac{\text{Res}(q,a_{j})}{a_{j}^{k+1}}\right).
\]

A perturbative calculation of error, assuming $\epsilon$ so small
that $\Delta a_{i}=\hat{a}_{i}-a_{i}$ satisfies $|\Delta a_{i}|\ll|a_{j}-a_{k}|$
for any $i,j,k$, follows. The perturbative calculation is based on
\[
\text{Res}(\hat{q},\hat{a}_{j})=\text{Res}(q,a_{j})\left(1-\sum_{i\neq j}\frac{\Delta a_{j}-\Delta a_{i}}{a_{j}-a_{i}}\right)+\mathcal{O}(\Delta a^{2})
\]
and
\[
\frac{1}{\hat{a}_{j}^{k+1}}=\frac{1}{a_{j}^{k+1}+ka_{j}^{k}\Delta a_{j}+\mathcal{O}(\Delta a_{j}^{2})}=\frac{1}{a_{j}^{k+1}}-\frac{\left(k+1\right)a_{j}^{k}\Delta a_{j}}{a_{j}^{2\left(k+1\right)}}+\mathcal{O}(\Delta a_{j}^{2}).
\]
These complete the first order perturbative calculation by implying
\begin{equation}
(q-\hat{q})_{k}=-\sum_{j=1}^{n}\frac{\text{Res}(q,a_{j})}{a_{j}^{k+1}}\left(\frac{(k+1)\Delta a_{j}}{a_{j}}+\sum_{i\neq j}\frac{\Delta a_{j}-\Delta a_{i}}{a_{j}-a_{i}}+\mathcal{O}(\Delta a^{2})\right).\label{eq:secn2.3-firstorder}
\end{equation}
Turning to condition number of roots of $p(z)=0$, we get the asymptotic
bound
\begin{equation}
|(q-\hat{q})_{k}|\lesssim\epsilon\sum_{j=1}^{n}\left|\frac{\text{Res}(q,a_{j})}{a_{j}^{k+1}}\right|\left(\frac{(k+1)\kappa(a_{j},p)}{|a_{j}|}+\sum_{i\neq j}\frac{\kappa(a_{i},p)+\kappa(a_{j,}p)}{|a_{j}-a_{i}|}\right).\label{eq:secn2.3-kappa}
\end{equation}
This bound suggests that the error in the $k$-th coefficient is dominated
by the root closest to $0$ in the limit $k\rightarrow\infty$. In
the transient phase, it suggests that some of the exterior roots may
dominate the error if they are sufficiently ill-conditioned. The latter
suggestion is not well-founded for a reason that will be presently
explained.

\begin{figure}
\begin{centering}
\subfloat[]{

\centering{}\includegraphics[scale=0.35]{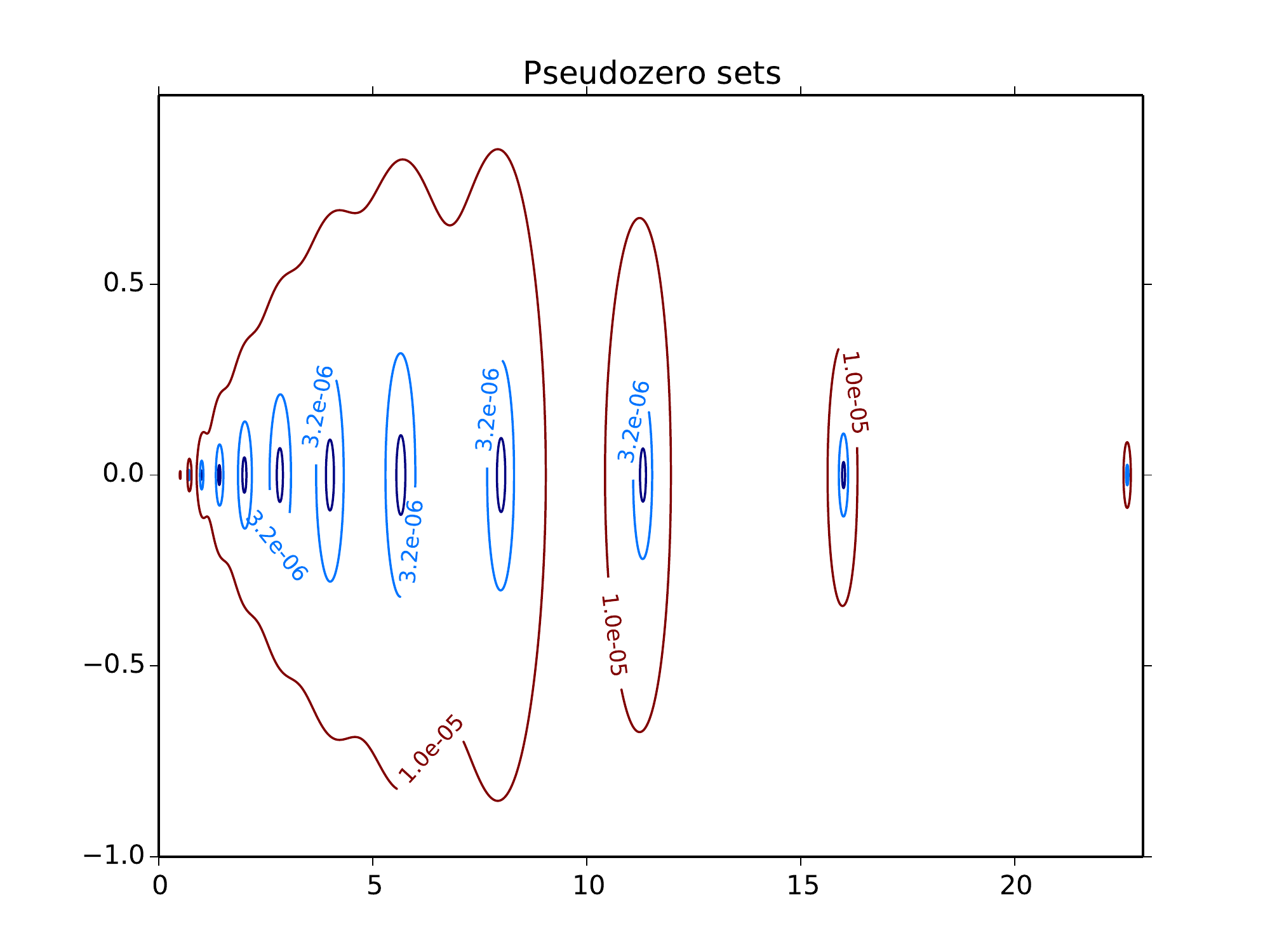}}\subfloat[]{\centering{}\includegraphics[scale=0.35]{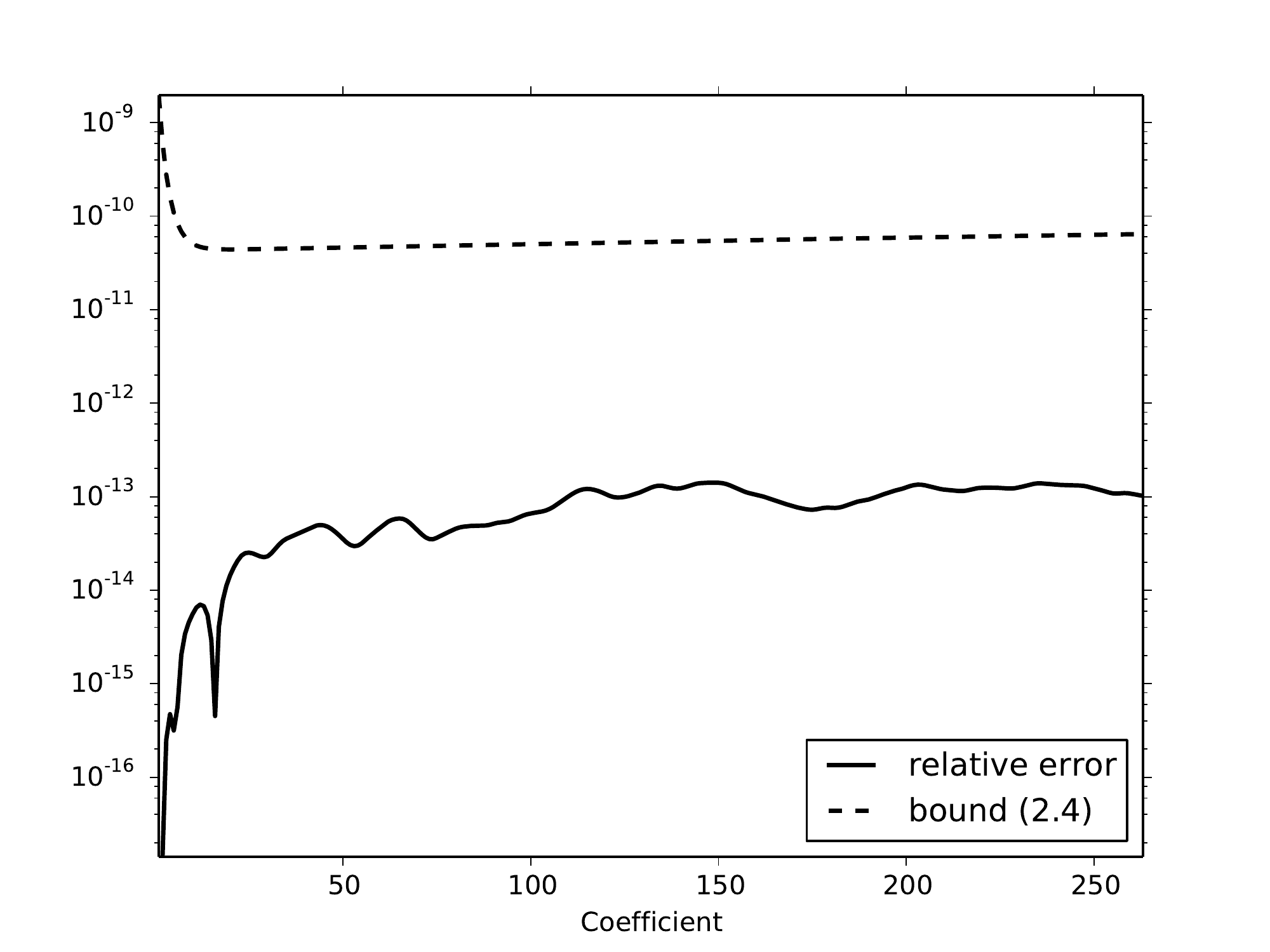}}
\par\end{centering}

\begin{centering}
\subfloat[]{\centering{}\includegraphics[scale=0.35]{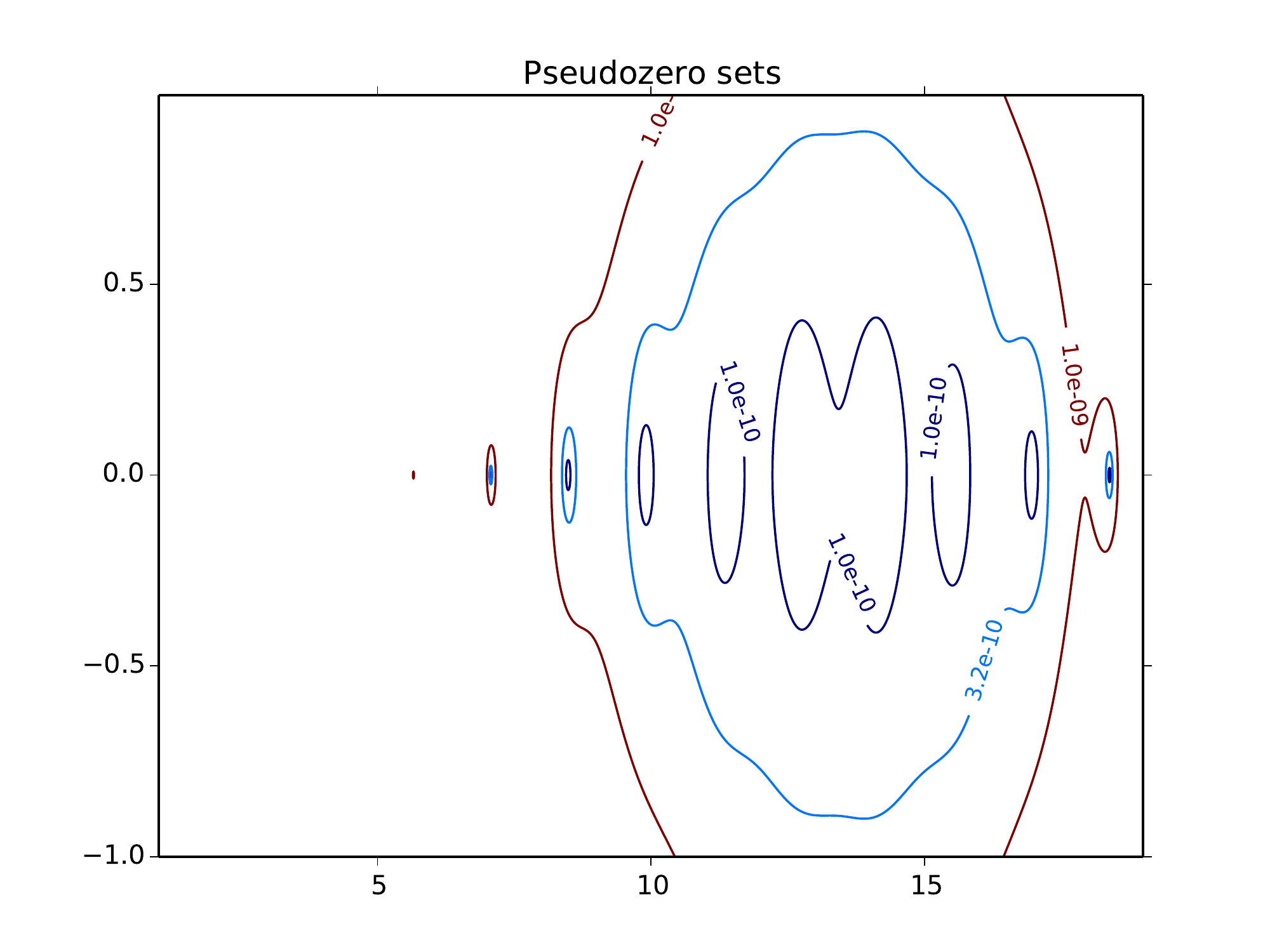}}\subfloat[]{

\begin{centering}
\includegraphics[scale=0.35]{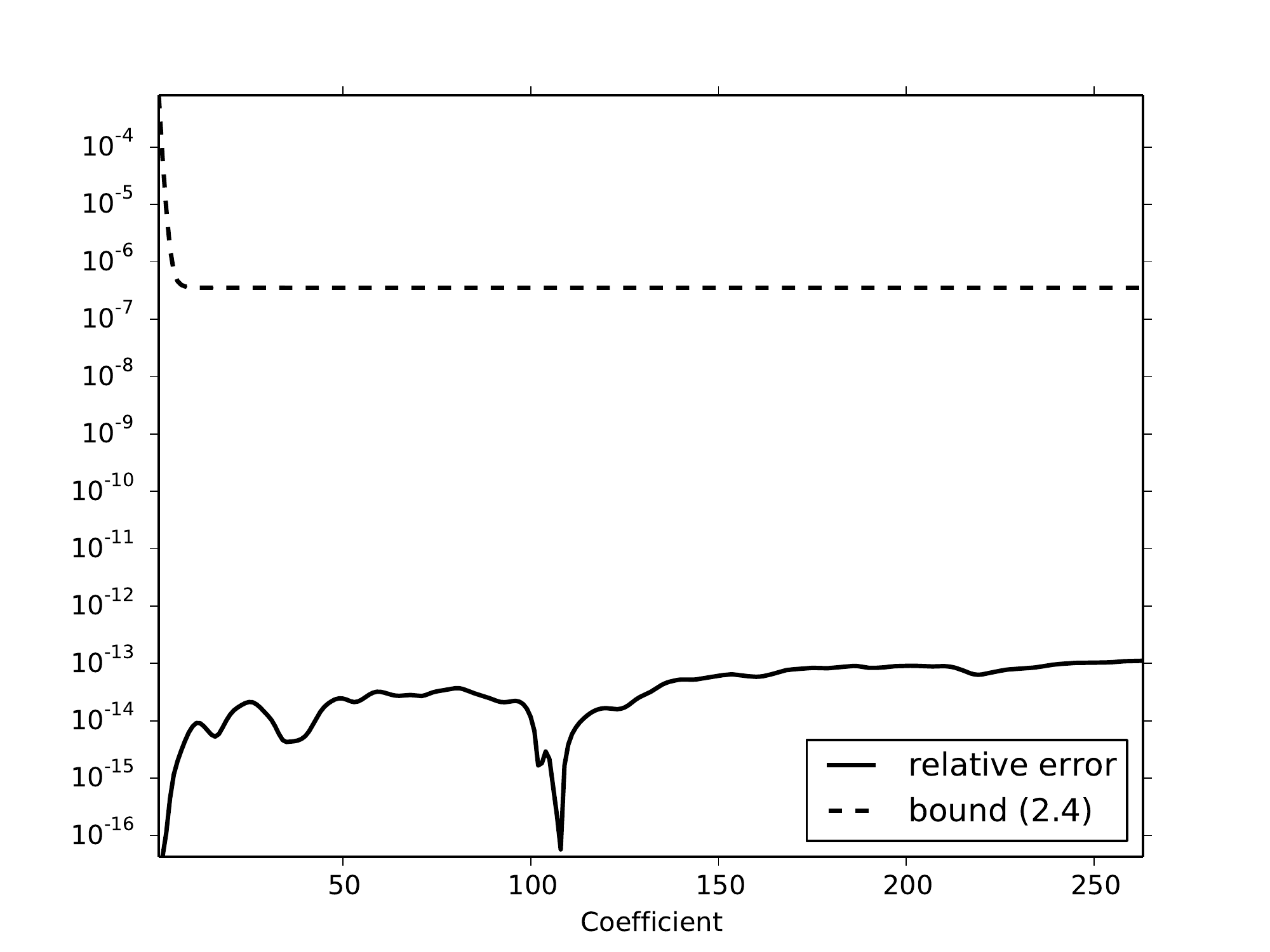}
\par\end{centering}

}
\par\end{centering}

\caption{Pseudozero sets and plots of relative error vs. coefficient for a
well-conditioned polynomial ((a) and (b)) and an ill-conditioned one
((c) and (d)), both of degree $13$. The bound (\ref{eq:secn2.3-kappa})
bounds absolute errors. It is converted to a bound on relative errors
in the plots. \label{fig:secn2.3-wilk-poly-inv} }
\end{figure}

Figure \ref{fig:secn2.3-wilk-poly-inv} compares the bound (\ref{eq:secn2.3-kappa})
(dashed line) to actual errors (solid line) for two examples. The
first example is $\prod_{i=-3}^{9}\left(x-2^{i/2}\right)$, implying
well-conditioned roots, and the second example is $\prod_{i=1}^{13}(x-i\sqrt{2})$,
implying ill-conditioned roots. Both examples are based on \cite{Wilkinson1984}.
In both examples, the bound (\ref{eq:secn2.3-kappa}) suggests transient
errors at the beginning which never materialize. The bound is highly
pessimistic for the ill-conditioned example. 

Part of the problem with the bound (\ref{eq:secn2.3-kappa}) is that
the condition numbers $\kappa(a_{j},p)$ can overestimate the perturbation
to the roots. But a more serious problem is that the errors $\Delta a_{i}$
in the first order error estimate (\ref{eq:secn2.3-firstorder}) are
highly correlated and this correlation is lost when they are bounded
using $\kappa(a_{j},p)$. Since $p$ and $\hat{p}$ are both monic
polynomials, the negative sums of their roots must equal $p_{n-1}$
and $\hat{p}_{n-1}$, respectively. Therefore no matter how large
each perturbation $\Delta a_{i}$ may be, their sum $\sum\Delta a_{i}$
must be of the order of machine precision, implying correlation between
the errors.

Such correlation between the errors $\Delta a_{i}$ is lost in the
asymptotic bound $|\Delta a_{i}|\lesssim\kappa(a_{i},p)$. Whether
the pseudozero plots contain information about correlations in the
errors is unknown. 

It is reasonable to expect a numerically stable algorithm for finding
roots to reproduce symmetric functions, such as the sum of all the
roots or the product of all the roots, accurately. However the algorithms
in current use progress from root to root, deflating the polynomial
every time a root is found. Perhaps for that reason they do not seem
to have this property. In particular algorithms that deflate using
the method of Theorem \ref{thm:secn2.1-fromcntoc0} will reproduce
the sum of the roots with accuracy but not the product of the roots.

\section{Error bounds and numerical stability}

The analysis given in this section uses techniques pioneered by Wilkinson
\cite{Wilkinson1961} and refined by Higham \cite{Higham1998,Higham2002}.
The application of the techniques is specialized to the inversion
of power series. Near the end of the section, we discuss the work
of Stewart \cite{Stewart1997} when comparing the errors that are
realized with the error bounds.

\subsection{Rounding error analysis}

To invert a power series as in 
\[
\frac{1}{1+b_{1}x+b_{2}x^{2}+\cdots}=1+c_{1}x+c_{2}x^{2}+\cdots
\]
the coefficients $c_{i}$ may be computed using 
\begin{eqnarray}
c_{1} & = & -b_{1}\nonumber \\
c_{2} & = & -b_{2}-c_{1}b_{1}\nonumber \\
 & \vdots\nonumber \\
c_{k} & = & -b_{k}-c_{1}b_{k-1}-\cdots-c_{k-1}b_{1}\label{eq:secn3.1-c1ck-recurrence}
\end{eqnarray}
The subtractions here are assumed to be left to right associative,
unlike Wilkinson's analysis of triangular back-substitution \cite{Wilkinson1961}
which assumes the opposite. Left to right associativity has the advantage
of preserving the Toeplitz structure of the matrices that arise in
error bounds. 

If we define $C_{n}$ and $T_{n}$ as 
\begin{equation}
C_{n}=\left(\begin{array}{c}
1\\
c_{1}\\
c_{2}\\
\vdots\\
c_{n}
\end{array}\right),T_{n=}\left(\begin{array}{ccccc}
1\\
b_{1} & 1\\
b_{2} & b_{1} & 1\\
\vdots & \vdots & \vdots & \ddots\\
b_{n} & b_{n-1} & b_{n-2} & \ldots & 1
\end{array}\right)\,\,\text{then}\,\, T_{n}^{-1}=\left(\begin{array}{ccccc}
1\\
c_{1} & 1\\
c_{2} & c_{1} & 1\\
\vdots & \vdots & \vdots & \ddots\\
c_{n} & c_{n-1} & c_{n-2} & \ldots & 1
\end{array}\right).\label{eq:secn3.1-CnTnTninv-defn}
\end{equation}
Here $T_{n}^{-1}$ is, like $T_{n}$, a Toeplitz matrix. In addition,
we have $T_{n}C_{n}={\bf e}_{1}$ where ${\bf e}_{1}$ is the vector
whose first component is $1$ and all others are $0$. In the recursion
(\ref{eq:secn3.1-c1ck-recurrence}) for computing $c_{k}$, the last
$c_{k-1}b_{1}$ term participates in only two arithmetic operations,
namely, the multiplication of $c_{k-1}$ and $b_{1}$ and the subtraction
of that product. Earlier terms participate in more subtractions and
the second term, which is $-c_{1}b_{k-1}$, participates in $k$ subtractions.
If the computed quantity is denoted $\hat{c}_{k}$, we may write
\[
\hat{c}_{k}=-b_{k}\left(1+\theta_{k+1}\right)-\hat{c}_{1}b_{k-1}\left(1+\theta_{k}\right)-\cdots-\hat{c}_{k-1}b_{1}\left(1+\theta_{2}\right).
\]
In other words, if $\hat{C}_{n}$ is the vector made up of $\hat{c}_{1},\ldots,\hat{c}_{n}$,
we have $\left(T_{n}+\Delta T_{n}\right)\hat{C}_{n}={\bf e}_{1}$
with $\bigl\vert\Delta T_{n}\bigr\vert\leq E_{n}$, where
\begin{equation}
E_{n}=\left(\begin{array}{ccccc}
0\\
\gamma_{2}|b_{1}| & 0\\
\gamma_{3}|b_{2}| & \gamma_{2}|b_{1}| & 0\\
\vdots & \vdots & \vdots & \ddots\\
\gamma_{n+1}|b_{n}| & \gamma_{n}|b_{n-1}| &  & \ldots & 0
\end{array}\right).\label{eq:secn3.1-En-defn}
\end{equation}
The identity 
\begin{equation}
\left(\hat{C}_{n}-C_{n}\right)=-T_{n}^{-1}\Delta T_{n}\left(\hat{C}_{n}-C_{n}\right)-T_{n}^{-1}\Delta T_{n}C_{n}\label{eq:secn3.1-chatn-c-identity}
\end{equation}
is the basis of the error bounds.

We may take norms of either side of (\ref{eq:secn3.1-chatn-c-identity})
and get 
\begin{equation}
|\hat{c}_{n}-c_{n}|\leq\bigl|\bigl|C_{n}-\hat{C}_{n}\bigr|\bigr|\leq\frac{\Bigl|\Bigl||T_{n}^{-1}|\: E_{n}\:|C_{n}|\Bigr|\Bigr|_{\infty}}{1-\Bigl|\Bigl||T_{n}^{-1}|\: E_{n}\Bigr|\Bigr|_{\infty}}.\label{eq:sec3.1-inf-norm-ineq}
\end{equation}
However, this bound is very poor. The coefficients of power series
are typically scaled badly, with terms increasing or decreasing at
a rapid rate. Norm-wise bounds are not of much use.

To get a component-wise bound, we go back to (\ref{eq:secn3.1-chatn-c-identity})
and take absolute values of both sides.
\begin{eqnarray*}
\left|\hat{C}_{n}-C_{n}\right| & \leq & \left|T_{n}^{-1}\right|\, E_{n}\,\left|\hat{C}_{n}-C_{n}\right|+\left|T_{n}^{-1}C_{n}\right|\, E\\
\left(I-\left|T_{n}^{-1}\right|\, E_{n}\right)\left|\hat{C}_{n}-C_{n}\right| & \leq & \left|T_{n}^{-1}C_{n}\right|\, E_{n}.
\end{eqnarray*}
Noting that the matrix $\left(I-\left|T_{n}^{-1}\right|\, E_{n}\right)$
is lower triangular with a non-negative inverse, we have the following
theorem.
\begin{thm}
If a power series is inverted using the recurrence (\ref{eq:secn3.1-c1ck-recurrence})
and left to right associativity, we have the error bound 
\begin{equation}
\left|\hat{C}_{n}-C_{n}\right|\leq\left(I-\left|T_{n}^{-1}\right|\, E_{n}\right)^{-1}\left|T_{n}^{-1}C_{n}\right|\, E_{n}.\label{eq:sec3.1-thm-3.1-ineq}
\end{equation}
\label{thm:secn3.1-cnhat-cn-bound}
\end{thm}

\subsection{Condition analysis and numerical stability}

If $p$ is a power series, $|p|$ denotes the power series with coefficients
replaced by their absolute values. Let $p$ and $q$ be power series
with constant terms equal to $1$ and
\[
pq=1.
\]
If $p$ is perturbed to $p+\Delta p$, where the constant term of
$\Delta p$ is $0$, suppose that $q$ gets perturbed to $q+\Delta q$.
We have 
\[
(p+\Delta p)(q+\Delta q)=1.
\]
It follows that 
\begin{eqnarray*}
p\Delta q & = & -q\Delta p-\Delta p\Delta q\\
\Delta q & = & -q^{2}\Delta p-q\Delta p\Delta q\\
|\Delta q| & \leq & |q^{2}\Delta p|+|q\Delta p|\,|\Delta q|\\
\left(1-|q|\Delta p|\right)|\Delta q| & \leq & |q^{2}\Delta p|
\end{eqnarray*}
All the coefficients of the power series $1/(1-|q\Delta p|)$ are
positive. Therefore we may multiply by that power series to get the
bound
\begin{equation}
|\Delta q|\leq\frac{|q^{2}\Delta p|}{1-|q\Delta p|}\leq\frac{|q^{2}|\,|\Delta p|}{1-|q|\,|\Delta p|}.\label{eq:secn3.2-cond-bound}
\end{equation}
We may take $|\Delta p$| to be 
\begin{equation}
\sum_{j=1}^{\infty}u|p_{j}|x^{j},\label{eq:secn3.2-Delta-p}
\end{equation}
where $u$ is the unit round-off, to obtain a bound on each entry
of $q$ using (\ref{eq:secn3.2-cond-bound}). Here it is significant
that the constant term of $\Delta p$ is zero. The conditioning bound
(\ref{eq:secn3.2-cond-bound}), with $|\Delta p|$ given by (\ref{eq:secn3.2-Delta-p}),
is sharp up to first order for each coefficient of $\Delta q$ with
a suitable choice of the signs of the coefficients of $|\Delta p|$.

Armed with this conditioning bound, we may consider the numerical
stability of the inversion of power series using the recurrence (\ref{eq:secn3.1-c1ck-recurrence}).
Theorem \ref{thm:secn3.1-cnhat-cn-bound} states that
\[
\left|\hat{C}_{n}-C_{n}\right|\leq\left(I-\left|T_{n}^{-1}\right|\, E_{n}\right)^{-1}\left|T_{n}^{-1}C_{n}\right|\, E_{n}.
\]
From the definitions of $C_{n}$ and $T_{n}^{-1}$ in (\ref{eq:secn3.1-CnTnTninv-defn})
as well as that of $E_{n}$ in (\ref{eq:secn3.1-En-defn}), we get
\[
\bigl|C_{n}-\hat{C}_{n}\bigr|\leq\frac{2(n+1)|q^{2}|\,|\Delta p|}{1-2(n+1)|q||\Delta p|}.
\]
Here we have used $\gamma_{k}<\gamma_{n+1}$ for $k\leq n$ and $\gamma_{n+1}\leq2(n+1)u$,
which assumes $(n+1)u<1/2$. This bound differs from the conditioning
bound (\ref{eq:secn3.2-cond-bound}) for each coefficient by only
a polynomial factor in $n$. Therefore inversion of power series using
back substitution is numerically stable.

\subsection{Numerical examples}

\begin{figure}
\subfloat[]{

\centering{}\includegraphics[scale=0.35]{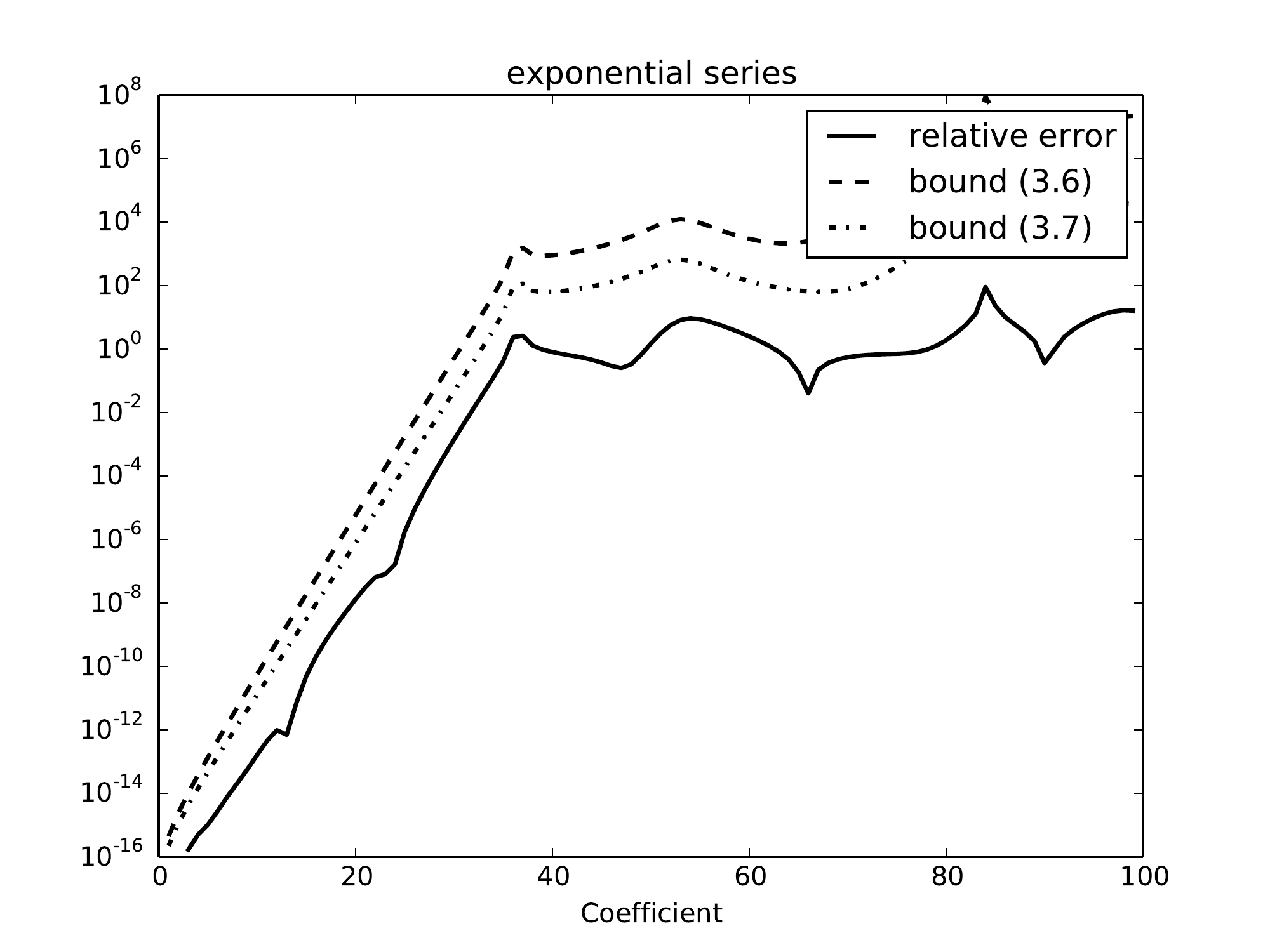}}\subfloat[]{\centering{}\includegraphics[scale=0.35]{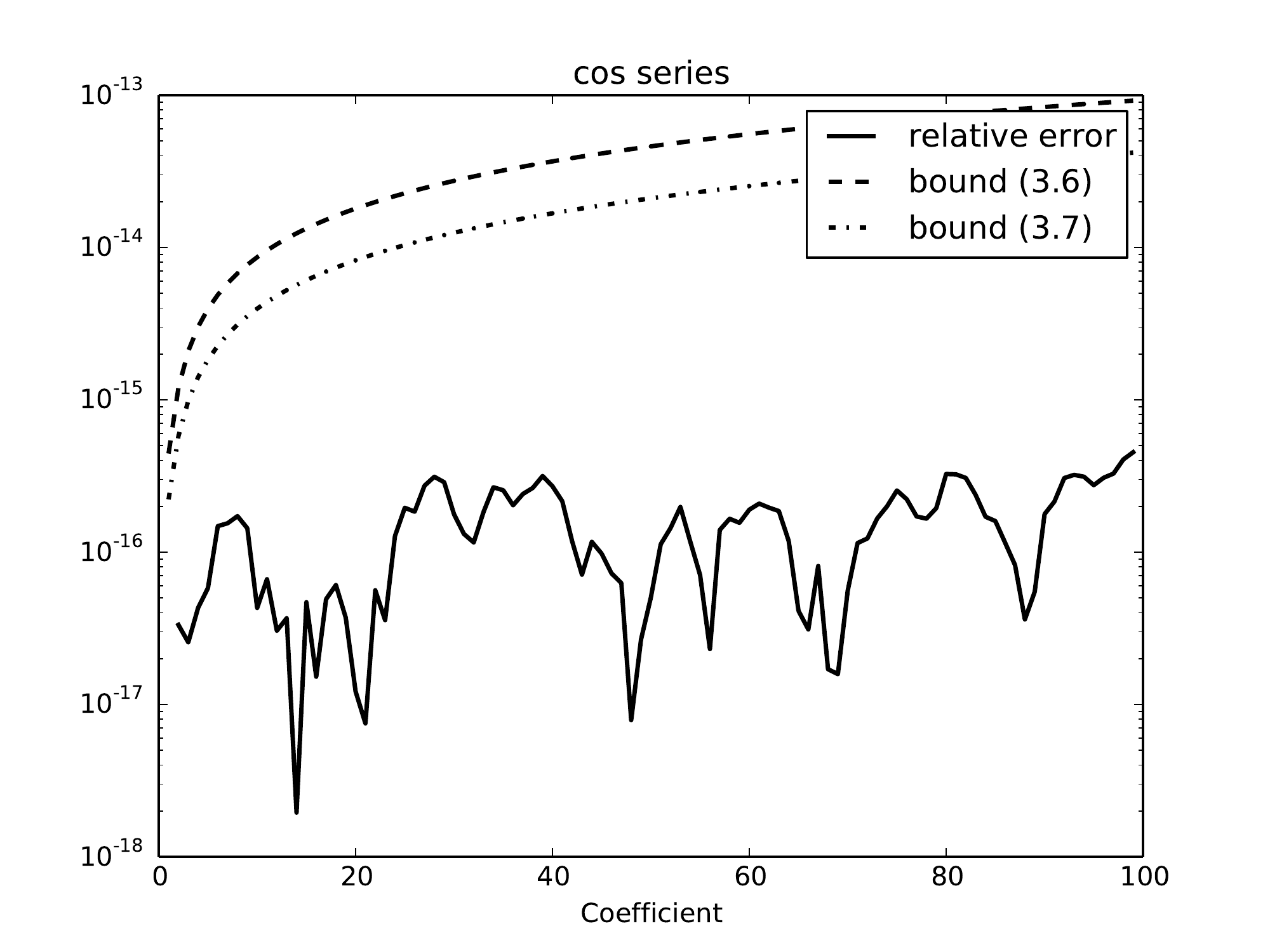}}

\subfloat[]{

\centering{}\includegraphics[scale=0.35]{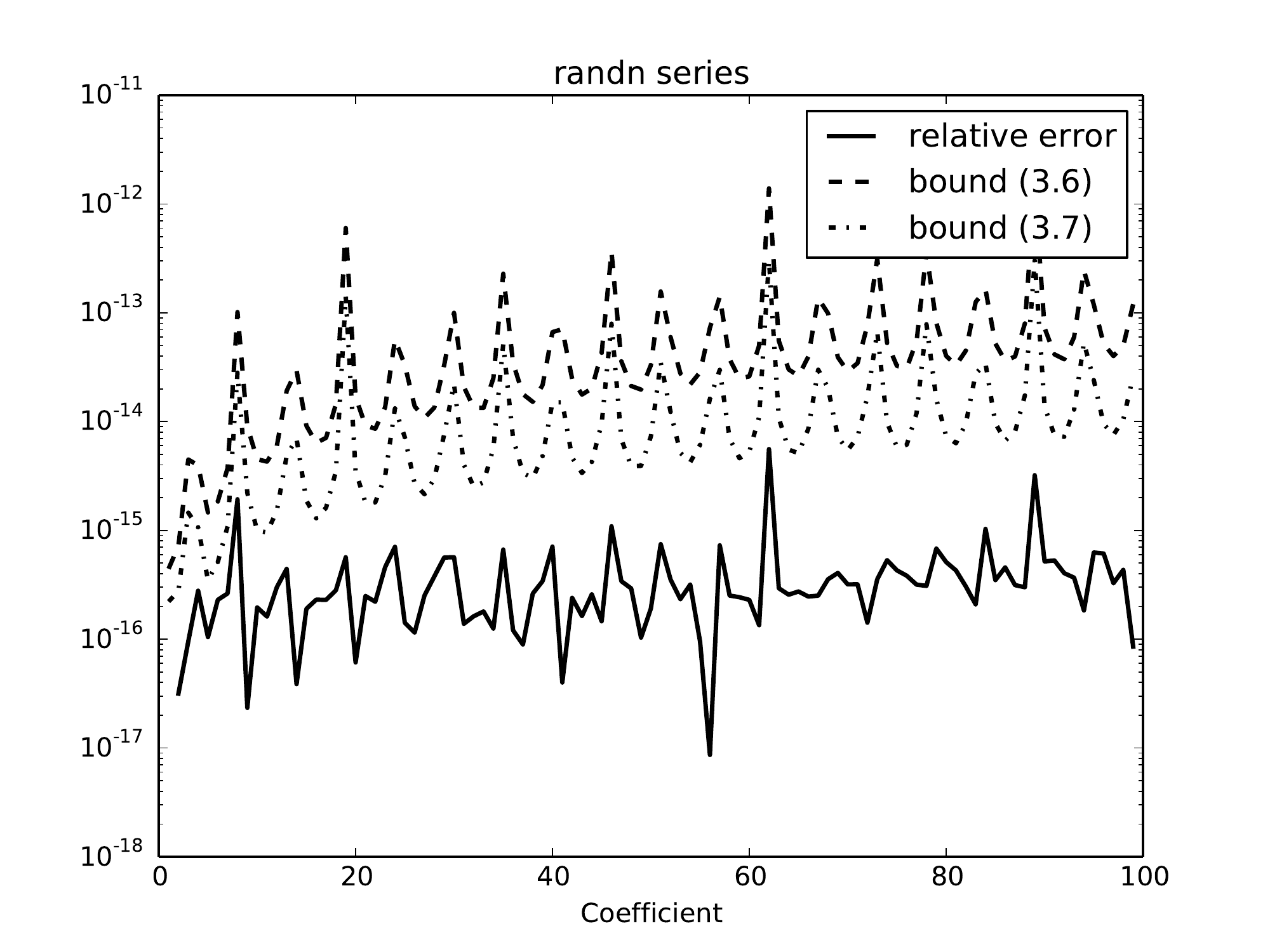}}\subfloat[]{\centering{}\includegraphics[scale=0.35]{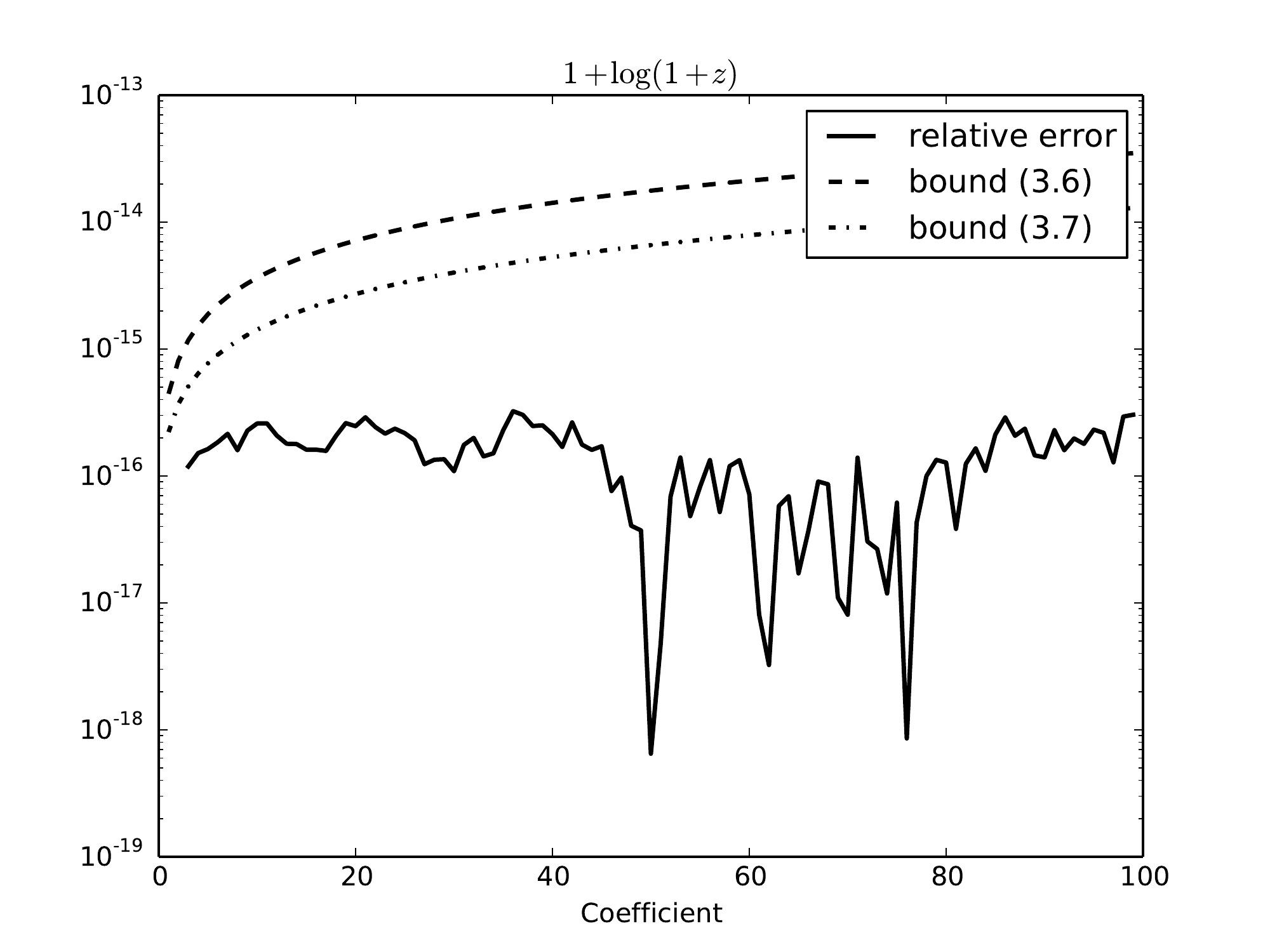}}

\caption{Rounding error bounds and actual rounding errors for four examples.
The bound of Theorem \ref{thm:secn3.1-cnhat-cn-bound} is on absolute
error. That bound is converted to a bound on relative error for bound
1 in each of the plots. Likewise, the bound of (\ref{eq:secn3.2-cond-bound}),
with $|\Delta p|$ given by (\ref{eq:secn3.2-Delta-p}), is converted
to bound 2. Each plot graphs relative error in the $n$-th coefficient
vs. $n$.\label{fig:secn3.3-numerical-examples}}
\end{figure}

Figure (\ref{fig:secn3.3-numerical-examples}) shows that the bounds
of Section 3.1 and 3.2 do quite well on four different examples. The
bounds themselves were computed using extended precision of $100$
digits. The actual relative error was computed by comparing the double
precision answers with extended precision answers. For inversion of
cosine, $1/\cos x$ in Figure (\ref{fig:secn3.3-numerical-examples})b,
the odd terms were ignored. It may be noted that the inverse cosine
series is one of the ways of defining Euler numbers. In the ``randn''
series, each $b_{i}$ in $p(x)=1+\sum_{i=1}^{\infty}b_{i}x^{i}$ is
an independent standard normal variable. 

Error bounds for inversion of triangular matrices are similar to that
of Theorem \ref{thm:secn3.1-cnhat-cn-bound}. However, they often
overestimate the error greatly \cite{Wilkinson1961}. In particular,
for many triangular matrices the relative error in the inverse appears
independent of the condition number. Here we discuss the work of Stewart
\cite{Stewart1997} and connect it to the inversion of power series. 

Consider the upper triangular matrix 
\begin{equation}
\left(\begin{array}{cc}
R & r\\
0 & \delta
\end{array}\right).\label{eq:secn3.3-stewart}
\end{equation}
If $\sigma$ is its smallest singular value, suppose $\sigma\geq\beta\delta$,
where $\beta\in[0,1]$ must hold. If $\beta$ is not too tiny, the
matrix is rank-revealing in the sense of Stewart. The last row of
this matrix may be rescaled to get 
\[
\left(\begin{array}{cc}
R & r\\
0 & 1
\end{array}\right),
\]
whose least singular value is denoted $\hat{\rho}$. If the least
singular value of $R$ is $\rho$, Stewart \cite{Stewart1997} has
proved that 
\[
\hat{\rho}\geq\frac{\beta\rho}{\sqrt{\beta^{2}+\rho^{2}}}.
\]
This bound may be interpreted as follows. If the matrix (\ref{eq:secn3.3-stewart})
is rank-revealing with a $\beta$ that is not too tiny, any significant
fall in the least singular value when we move from $R$ to that matrix
must be due to the smallness of $\delta$. The smallness of $\delta$
can be easily eliminated by rescaling the last row to get a matrix
whose condition number $\hat{\rho}$ is only moderately smaller than
$\rho$ the condition number of $R$. On the other hand, if the best
possible $\beta$ is quite tiny, it may mean that the ill-conditioning
of the matrix (\ref{eq:secn3.3-stewart}) is hidden within the correlations
between rows in a way that may not be eliminated so easily. If each
one of the principal submatrices of a matrix is rank revealing, any
ill-conditioning is almost entirely removed by rescaling rows explaining
Wilkinson's observation.

Many triangular matrices are not rank-revealing. For example, random
triangular matrices are not rank-revealing with probability $1$ as
proved in \cite{ViswanathTrefethen1998}. However, Stewart \cite{Stewart1997}
argues intuitively that the triangular matrices that arise in Gaussian
elimination and QR factorization are likely to be rank revealing.
His argument is that if a matrix is rank deficient, Gaussian elimination
and QR will break down with a $0$ on the diagonal. If it is nearly
rank deficient, continuity suggests that a very small entry must appear
on the diagonal indicating its rank deficiency. Pivoting makes either
factorization more apt to be rank revealing.

To connect Stewart's analysis to power series, we shall assume that
$p(x)=1+\sum b_{i}x^{i}$ has a radius of convergence $R$ equal to
$1$. Any finite radius of convergence can be turned into $1$ by
the change of variables $x\leftarrow x/R$. Assuming $R=1$, the matrix
$T_{n}$ of (\ref{eq:secn3.1-CnTnTninv-defn}) is rank-revealing if
and only if its least singular value is $\mathcal{O}(1)$. The least
singular value of $T_{n}$ is $\mathcal{O}\left(1\right)$ if and
only if the greatest singular value of $T_{n}^{-1}$is $\mathcal{O}(1)$,
which is true if and only if the entries $c_{i}$ of $T_{n}^{-1}$
in (\ref{eq:secn3.1-CnTnTninv-defn}) are $\mathcal{O}(1)$. Since
$c_{i}$ are the coefficients of the power series of $1/p(x)$, we
have that $T_{n}$ is rank revealing in the sense of Stewart if and
only if the radius of convergence of $1/p(x)$ is $1$ or greater. 

If the equation $p(z)=0$ has a solution with $|z|<1$ in the complex
plane, the matrix $T_{n}$ will not be rank-revealing. The example
of Figure \ref{fig:secn3.3-numerical-examples}d, $p(z)=1+\log(1+z)$
has a zero at $z=1-1/e$ and the corresponding matrix $T_{n}$ is
not rank-revealing. If in fact the radius of converge of $p(z)$ is
$1$ and there is no zero with $|z|<1$, the matrix $T_{n}$ will
be rank-revealing but its condition number will be $\mathcal{O}(1)$.
Within the scope of the analysis given by Stewart, the situation where
the actual relative errors are much smaller than the conditioning
bound appears unlikely. The good agreement between the bounds and
the actual errors in Figure \ref{fig:secn3.3-numerical-examples}
is the rule rather than the exception.

\section{Conclusions}

In this article, we have considered the inversion of power series
with particular attention to the special case of inverting polynomials.
Essential background is provided by the classic work of Wilkinson
\cite{Wilkinson1961} on inversion of triangular systems. 

We found and explicated a subtle numerical instability that arises
when factors corresponding to known roots are deflated from polynomials.
This instability has occurred in the computation of spectral differentiation
matrices. The suggestion that polynomial root finding algorithms such
as Jenkins-Traub may be more accurate without the deflation step merits
further investigation. 

The rounding error analysis and the condition analysis of power series
inversion imply numerical stability. In addition, the error bounds
that result from the analysis are not unduly pessimistic, as happens
for certain other triangular systems.

\section{Acknowledgements}

This research was partially supported by NSF grants DMS-1115277 and
SCREMS-1026317.

\bibliographystyle{plain}
\bibliography{references}

\end{document}